\documentstyle[art10,mssymb]{article}

\newcommand{\arrow}{\longrightarrow}
\newcommand{\Z}{{\Bbb Z}}
\newcommand{\C}{{\Bbb C}}
\newcommand{\R}{{\Bbb R}}

\newcommand{\g}{{\goth g}}

\newcommand{\restrict}[1]{{|_{{\phantom{|}\!\!}_{#1}}}}
\renewcommand{\c}[1]{{\cal #1}}

\renewcommand{\phi}{\varphi}
\renewcommand{\epsilon}{\varepsilon}
\renewcommand{\geq}{\geqslant}
\renewcommand{\leq}{\leqslant}

\newcommand{\comment}[1]{{}}

\def\blacksquare{\hbox{\vrule width 4pt height 4pt depth 0pt}}

\makeatletter

\newcounter{lemma}[section]
\renewcommand{\thelemma}{{Lemma \thesection.\arabic{lemma}}}
\newcommand{\lemma}{%
     \refstepcounter{lemma}
     {\bf \thelemma:\ }}

\newcounter{claim}[section]
\renewcommand{\theclaim}{{Claim \thesection.\arabic{claim}}}
\newcommand{\claim}{%
     \refstepcounter{claim}
     {\bf \theclaim:\ }}

\newcounter{sublemma}[section]
\newcounter{corollary}[section]
\renewcommand{\thecorollary}{{Corollary \thesection.\arabic{corollary}}}
\newcommand{\corollary}{%
     \refstepcounter{corollary}
     {\bf \thecorollary:\ }}

\newcounter{theorem}[section]
\renewcommand{\thetheorem}{{Theorem \thesection.\arabic{theorem}}}
\newcommand{\theorem}{%
     \refstepcounter{theorem}
     {\bf \thetheorem:\ }}

\newcounter{proposition}[section]
\newcounter{definition}[section]
\renewcommand{\thedefinition}
       {{Definition \thesection.\arabic{definition}}}
\newcommand{\definition}{%
     \refstepcounter{definition}
     {\bf \thedefinition:\ }}

\newcounter{example}[section]
\newcounter{problem}[section]
\newcounter{question}[section]
\newcommand{\eqref}[1]{(\ref{#1})}

\newcommand{\ps@verbit}{%
  \renewcommand{\@oddhead}{%
          \scriptsize
          {Cohomology of compact hyperk\"ahler... (applications)}
          \hfil\tiny {August 95}}
  \renewcommand{\@evenhead}{\@oddhead}
  \renewcommand{\@oddfoot}{\hfil\thepage\hfil}
  \renewcommand{\@evenfoot}{\@oddfoot}}

\makeatother

\begin{document}

\begin{center}
\Large\bf
Cohomology of compact hyperk\"ahler manifolds\\
and its applications.
\end{center}

\centerline{Mikhail Verbitsky,}
\centerline{verbit@math.harvard.edu}

\hfill

\hfill

\hfill

\hfill

\hspace{0.1\linewidth}\begin{minipage}[t]{0.8\linewidth}
{\bf Abstract}. This article contains
a compression of results from \cite{_main_},
with most proofs omitted.
We prove that every two points of the connected
moduli space of holomorphically symplectic manifolds
can be connected with so-called ``twistor lines'' --
projective lines holomorphically embedded to the moduli space
and corresponding to the hyperk\"ahler structures.
This
has interesting implications for the geometry of compact
hyperk\"ahler manifolds and of holomorphic vector
bundles over such manifolds.
\end{minipage}

\hfill


\section{Lie algebra action.}


We refer to \cite{_main_} for details of definitions
and missing proofs.
A hyperk\"ahler manifold is a Riemannian manifold
$M$ equipped with three complex structures $I$,
$J$ and $K$, such that $I\circ J=-J\circ I=K$ and $M$
is K\"ahler with respect to $I$,
$J$ and $K$. Relations between $I$, $J$ and $K$
imply that there is is an action of quaternions in its tangent
space. Consequently, there is a multiplicative action of
$SU(2)$ on the algebra of differential forms. This action
commutes with Laplacian. Hence there is a canonical
action of $SU(2)$ on cohomology of $M$.

Let $M$ be a complex manifold which admits a hyperk\"ahler
structure. A simple linear-algebraic argument implies that
$M$ is equipped with a holomorphic symplectic form. Calabi-Yau theorem
shows that, conversely, every compact holomorphically symplectic
K\"ahler manifold admits a hyperk\"ahler structure, which is uniquely
defined by these data. Further on, we consider only holomorphically symplectic
manifolds which are compact and of K\"ahler type.
For simplisity of statements, we assume also that
\[ \dim H^{2,0}(M)=1, \mbox{\ \  and \ \ }H^1(M) =0, \]though these
assumptions are not necessarily for most results.

\hfill

The algebraic structure on $H^*(M)$ is studied using the
general theory of Lefschetz-Frobenius algebras, introduced
in \cite{_Lunts-Loo_}.

Let $A= \bigoplus\limits^{2d}_{i=0} A_i$
be a graded commutative associative algebra over a field of
characteristic zero. Let $H\in End(A)$ be a linear endomorphism
of $A$ such that for all $\eta \in A_i$, $H(\eta)= (i-d) \eta$.

For all $a\in A_2$, denote by $L_a:\;\; A\arrow A$ the linear map which
associates with $x\in A$ the element $ax\in A$.
The triple $(L_a, H, \Lambda_a)\in End(A)$ is called {\bf a Lefschetz triple}
if

\[ [ L_a, \Lambda_a] = H,\;\; [ H, L_a ] = 2 L_a, \;\;
   [ H, \Lambda_a] = -2 \Lambda_a.
\]
A Lefschetz triple establishes a representation of the
Lie algebra $\goth{sl}(2)$ in the space $A$. For cohomology
algebras, this representation
arises as a part of Lefschetz theory.
In a Lefschetz triple, the endomorphism $\Lambda_a$ is
uniquely defined by the element $a\in A_2$
(\cite{_Bou:Lie_} VIII \S 3). For arbitrary
$a\in A_2$, $a$ is called {\bf of Lefschetz type} if the
Lefschetz triple $(L_a, H, \Lambda_a)$ exists. If $A= H^*(X)$
where $X$ is a compact complex manifold of K\"ahler type,
then all K\"ahler classes $\omega\in H^2(M)$ are elements
of Lefschetz type. As one can easily check,
the set $S\subset A_2$ of all elements of Lefschetz type
is Zariski open in $A_2$.

\hfill

\definition 
A Lefschetz-Frobenius algebra is a
Frobenius graded supercommutative algebra which admits a
Lefschetz triple.

\hfill

\definition 
Let $A$ be a Lefschetz-Frobenius algebra.
The structure Lie algebra $\g(A)\subset End(A)$ is a Lie subalgebra
of $End(A)$ generated by $L_a$, $\Lambda_a$, for all
elements of Lefschetz type $a\in S$.

\hfill

Let $M$ be a compact hyperk\"ahler manifold with the complex
structures $I$, $J$, $K$. Consider the K\"ahler forms
$\omega_I$, $\omega_J$, $\omega_K$ associated with these
complex structures. Let $\rho_I:\; \goth{sl}(2)\arrow End(H^*(M))$,
$\rho_J:\; \goth{sl}(2)\arrow End(H^*(M))$,
$\rho_K:\; \goth{sl}(2)\arrow End(H^*(M))$
be the corresponding Lefschetz homomorphisms.
Let $\goth a\subset End(H^*(M))$ be the minimal
Lie subalgebra which contains images of $\rho_I$,
$\rho_J$, $\rho_K$. The algebra $\goth a$ was computed
explicitely in \cite{_so5_on_cohomo_}.

\hfill

\theorem \label{_so_5_Theorem_} 
(\cite{_so5_on_cohomo_})
The Lie algebra $\goth a$ is naturally isomorphic to $\goth{so}(4,1)$.

\hfill

This statement can be regarded as a ``hyperk\"ahler Lefschetz theorem''.
Indeed, its proof parallels the proof of Lefschetz theorem.

\hfill


Using \ref{_so_5_Theorem_}, we compute the structure Lie algebra
of $H^*(M)$.

\hfill

\theorem \label{_structure_alge_for_coho_hyperkahe_Theorem_} 
(\cite{_main_}, Theorem 11.1)
Let $M$ be a compact holomorphically symplectic manifold. Assume
that $\dim H^{2,0}(M)$
Let $n= \dim(H^2(M))$. Let $\g(A)$ be a structure Lie algebra
for $A=H^*(M)$. Then $\g(A)$ is isomorphic to
$\goth{so}(4, n-2)$.\footnote{This isomorphism can be made canonical.
The Lie algebra $\g(A)$ is isomorphic to $\goth{so}(V\oplus \c H)$
where $V$ is the linear space $H^2(M, \R)$ equipped with the
natural pairing of a signature $(3,n-3)$ (\cite{_Beauville_}
Remarques, p. 775; see also \ref{_H-R_form_defi_Theorem_}), and $\c H$
is 2-dimensional vector space with hyperbolic quadratic form.}

\hfill

Let $H^*_r(M)$ be a sub-algebra of $H^*(M)$ generated by $H^2(M)$.
It is easy to see that $\g(A)$ acts on $H^*_r(M)$, and $H^*_r(M)$
is an irreducible representation of $\g(A)$. Moreover, multiplicative
structure in $H^*_r(M)$ is easily recovered from an action of $\g(A)$.
Using the general knowledge of representations of $\goth{so}(n)$,
we obtain exact knowledge of the multiplicative structure
of $H^*_r(M)$. In particular, we obtain the
following theorem.

\hfill

\theorem \label{_S^*H^2_is_H^*M_intro-Theorem_} 
(\cite{_main_}, Theorem 15.2)
Let $\dim_\C M=2n$.
Then

\[\bigg\{\begin{array}{lr}
 H^{2i}_r(M)\cong S^i H^2(M)&
 \mbox{\ \ for $i\leq n$, and}\\
 H^{2i}_r(M)\cong S^{2n-i} H^2(M) &
 \mbox{\ \ for $i\geq n$ }
 \end{array}
\]

\section{The Riemann-Hodge pairing.}


Let $M$ be a compact holomorphically symplectic manifold
f K\"ahler type, satisfying
\[ \dim H^{2,0}(M)=1. \]
In \cite{_Beauville_} Remarques, p. 775, Beauville introduces
canonical 2-form on $H^2(M)$, of signature
$(n-3,3)$, where $n=\dim H^2(M)$. Throughout the paper \cite{_main_},
this form was called {\bf the Riemann-Hodge
pairing}.\footnote{More accurately, this form should be
called Bogomolov-Beauville form.
The author was unaware of Beauville's
remark, and did not understand the part of Bogomolov's
paper (\cite{_Bogomolov_}) where this form is also introduced.}
In \cite{_main_}, this form was described via the action of
$SU(2)$ on $H^2(M)$.

\hfill

 Let $\omega$ be a K\"ahler class on $M$
such that \[\int_M \omega^{\dim_\C M}=1,\] and $(I,J, K, (\cdot,\cdot))$
be the corresponding hyperk\"ahler structure. Let
\[ (\cdot,\cdot)_{Her}:\; H^2(M,\C)\times H^2(M,\C)\arrow \C\]
be a positively Hermitian form on second cohomology of $M$ which
corresponds to the Riemannian structure $(\cdot,\cdot)$.
Let $H^2(M) = H^{inv}(M)\oplus H^{+}(M)$ be a decomposition
such that
$H^{inv}(M)$ consists of all $SU(2)$-invariant 2-forms,
and $H^{+}(M)$ is the complementary $SU(2)$-invariant subspace.
Let $(\cdot,\cdot)_{\c H}$ be the form which is equal to
$(\cdot,\cdot)_{Her}$ on $H^{+}(M)$ and $-(\cdot,\cdot)_{Her}$
on $H^{inv}(M)$

\hfill

\theorem \label{_H-R_form_defi_Theorem_} 
(\cite{_main_}, Theorem 6.1, cf. \cite{_Beauville_} Remarques, p. 775)
The form $(\cdot,\cdot)_{\c H}$ is independent
from the choice of the complex and K\"ahler structure on $M$.

\hfill

The form $(\cdot,\cdot)_{\c H}$ is used in the proof
of \ref{_structure_alge_for_coho_hyperkahe_Theorem_}.

\hfill

Let $\rho_I:\; \goth{u}(1)\arrow End(H^*(M))$ be a map for which
$z\in \goth{u}(1)$ acts on $H^{p,q}(M)$ by $(p-q)z$.
Clearly, the action of $\goth{u}(1)$ on $H^2(M)$ respects the form
$(\cdot,\cdot)_{\c H}$.
Let $\g_M\subset End(H^*(M))$ be a Lie algebra generated by the images
of $\rho_I$ for all complex structures $I$ on $M$.
Let $V$ denote the linear space $H^2(M)$ equipped with
bilinear form $(\cdot,\cdot)_{\c H}$.
By \ref{_H-R_form_defi_Theorem_}, the action of $\g_M$ on
$V$ preserves $(\cdot,\cdot)_{\c H}$.
This defines a Lie algebra homomorphism
$\Gamma:\;\g_M \arrow  \goth{so}(V)$. The following theorem
is a chief tool in proving the Mirror Conjecture
for a compact holomorphically symplectic manifold.

\hfill

\theorem \label{_embe_Mum_Tate_to_so_Theorem_} 
The map $\Gamma:\; \g_M \arrow  \goth{so}(V)$ is an isomorphism.

{\bf Proof:} \cite{_main_}, Theorem 13.1, 13.2. $\blacksquare$

\hfill

The Lie algebra $\g(A)\subset End(H^*(M))$ is equipped with a natural
grading, induced by the grading on $H^*(M)= \oplus H^i(M)$.
Let $k$ be the one-dimensional Lie subalgebra of $End(H^*(M))$
spanned by $Id$.

\hfill

\theorem \label{_Mum_Tate_is_g_0_Theorem_} 
(\cite{_main_}, Theorem 13.2)
The Lie subalgebra \[ \g_M\oplus k\subset End(H^*(M))\]
coinsides with the grading-zero part of $\g(A)$.

$\blacksquare$

\hfill

We have a {\bf period map}
$P_c:\; Comp \arrow {\Bbb P}H^2(M, \C)$ associating a line
$H^{2,0}_I(M) \subset H^2(M, \C)$ to a complex structure $I$.
Complexifying $H^2(M, \R)$, we can consider
$(\cdot,\cdot)_{\c H}$ as a complex-linear,
complex-valued form on $H^2(M, \R)$. For all $I\in Comp$,
$P_c(I)$ belongs to a conic hypersurface $C \subset {\Bbb P}H^2(M, \C)$,

\[ C = \{ l \;\; |\;\; (l,l)_{\c H} =0\}. \]
Torelli principle (proved by Bogomolov in the case of holomorphically
symplectic manifolds, \cite{_Bogomolov_})
implies that $P_c:\; Comp \arrow C$ is etale.

\hfill

Let $\underline{\c H}= \oplus H^{p,q}(M)$ be a variation of Hodge structures
(VHS) on Comp associated with the total cohomology space of $M$.
\ref{_embe_Mum_Tate_to_so_Theorem_} implies that there exist a VHS
$\c H$ on $C$, such that $\underline {\c H}$
is a pullback of a variation of Hodge structures $\c H$:
$\underline {\c H} = P^*_c(\c H)$.  Let $G_M$ be the Lie group
associated with $\g_M$,
$G_M = Spin\bigg(H^2(M, \R), (\cdot,\cdot)_{\c H}\bigg)$.
The set $C$ is equipped with a natural action of a group $G_M$.
This group also acts in the total cohomology space $H^*(M)$
of $M$. This defines an equivariant structure in the
bundle $\c H$. The chief idea used in the proof
of Mirror Symmetry is the following theorem:

\hfill

\theorem \label{_VHS_equi_Intro_Theorem_} 
The VHS $\c H$ is $G_M$-equivariant, under the natural
action of $G_M$ on $C$ and $\c H$.

{\bf Proof:} See \cite{_V:Mirror_}, Theorem 2.2.
$\;\;\blacksquare$

\hfill

To make this statement more explicit, we recall that the variation
of Hodge structures is a flat bundle, equipped with a real structure
and a holomorphic filtration (Hodge filtration),
which is complementary to its complex adjoint filtration.
Then, \ref{_VHS_equi_Intro_Theorem_} says
that the action of $G_M$ on $\c H$
maps flat sections to flat sections, and preserves the real structure
and the Hodge filtration.


\section{The twistor lines.}


The main technical tool used in the text of \cite{_main_} is
results about (coarse, marked) moduli space $Comp$ of complex structures
on a holomorphically symplectic manifold $M$.
Let $\omega$ be a K\"ahler class on $M$ and $\c H=(I,J,K, (\cdot,\cdot))$
be the corresponding hyperk\"ahler structure. Then,
for every triple of real numbers $(a,b,c), a^2+b^2+c^2=1$,
the operator $aI+bJ+cK$ defines an integrable complex
structure%
\footnote{This complex structure is
called {\bf a complex structure induced by a hyperk\"ahler
structure}.}
on $M$.
Identifying the set of such triples
with $\C P^1$, we obtain a map
$\C P^1\stackrel{i_{\c H}}{\hookrightarrow} Comp$
where $Comp$ is a connected component of the coarse
moduli space of $M$.
The following claim is easy.

\hfill

\claim 
The map $i_{\c H}$ is a
holomorphic embedding of complex analytic varieties.

\hfill

Let $P:\; Comp \arrow C$ be the period map, assigning
to a complex structure $I$ a line $H^{2,0}(M,I)$.
Let $C\subset {\Bbb P}^1(H^2(M,\C)= P(Comp)$.
 According to \cite{_Beauville_}, $P$ is etale.
The projective line $i_{\c H}(\C P^1)\subset Comp$ is called
{\bf a twistor line}, and is denoted by $R_{\c H}$.
The following theorem was, regrettably, omitted in \cite{_main_},
though all necessary tools were developed for its proof.
For conceptual understanding of our argument, this
theorem is indispensable.

\hfill

\theorem \label{_twistor_connect_Theorem_} 
Let $I_1, I_2\in Comp$. Then there exist a sequence of
intersecting twistor lines which connect $I_1$ with $I_2$.

{\bf Proof:}
To prove \ref{_twistor_connect_Theorem_}$'$, we have to show that
a set $\tilde{\c L_0}$ of all twistor lines $i_{\c H_0}(\C P^1)$
which are connected to $i_{\c H}(\C P^1)$ with
intersecting twistor lines is open.
Since $P:\; Comp \arrow C$ is etale, it suffices to show that
$I_1, I_2$ can be connected with twistor lines $l_i$
such that $P(l_i)$ intersect $P(l_{i+1})$.

With every twistor line $R_{\c H}$, we associate
a 3-dimensional plane $\ell_{\c H}\subset H^2(M,\R)$ which is
spanned by the K\"ahler classes $\omega_I$, $\omega_J$, $\omega_K$.
A linear algebraic argument shows that the twistor
lines $R_{\c H_1}$ and $R_{\c H_2}$ intersect if
and only if $\dim (\ell_{\c H_1}\cap \ell_{\c H_1})\geq 2$.
Hence we need to show that

\hfill

{\bf \ref{_twistor_connect_Theorem_}$'$} 
Each pair of twistor lines $R_{\c H}$, $R_{\c H'}$
can be connected with
a sequence of twistor lines $R_{\c H}=R_{\c H_1}$, ...,
$R_{\c H_n}=R_{\c H'}$ such that
$\dim (\ell_{\c H_i}\cap \ell_{\c H_{i+1}})\geq 2$.

\hfill

\lemma \label{_Kah_open_Lemma_} 
Let $H$ be a hyperk\"ahler structure on $M$, $i_{\c H}(\C P^1)\subset Comp$
be the set of all induced complex structures, and $Kah(\c H)$ be the set
of all K\"ahler classes corresponding to $L\in i_{\c H}(\C P^1)$. Then
$Kah(\c H)$ is open in $H^2(M, \R)$.

{\bf Proof:} \cite{_main_}, Claim 6.6 $\;\;\blacksquare$

\hfill

Let $\c L$ be the space of all triples
$\omega_I$, $\omega_J$, $\omega_K$ in $H^2(M)$ which are orthonormal
with respect to the pairing $(\cdot,\cdot)_{\c H}$ of
\ref{_H-R_form_defi_Theorem_}, and $Hyp$ be the connected
component of the set of all hyperk\"ahler structures.
Let $P_h:\; Hyp \arrow \c L$ be the natural period map.
Comparing dimensions and using Calabi-Yau, we observe
that $P_h$ is etale. Let $\c L_0$ be the
space of twistor lines corresponding to $\tilde{L_0}$.
Using \ref{_Kah_open_Lemma_}, we find that the differential
of $P_h\restrict{\c L_0}$ is surjective. Therefore,
$\c L_0$ is open in $\c L$, and $\tilde{\c L_0}$ is open
in the set of all twistor lines. This proves
\ref{_twistor_connect_Theorem_}. $\blacksquare$

\hfill


\section{An outline of proofs.}


Let $(I,J, K, (\cdot,\cdot))$ be a hyperk\"ahler structure
on $M$. One can check that the cohomology classes
$\omega_I$, $\omega_J$, $\omega_K\in H^2(M,\R)$ are orthogonal
with respect to the pairing $(\cdot,\cdot)$. Let $Hyp$ be the classifying
space of the hyperk\"ahler structures on $M$. Let
$P_{hyp}:Hyp\arrow H^2(M)\times H^2(M)\times H^2(M)$
be the map which associates with the hyperk\"ahler structure
$\c H= (I, J, K, (\cdot,\cdot))$ the triple
$(\omega_I,\omega_J,\omega_K)$. Then the image of $P_{hyp}$ in
$H^2(M)\times H^2(M)\times H^2(M)$ satisfies

\begin{equation}\label{_image_of_P_hyp_Equation_}
   \forall (x,y,z)\in im P_{hyp}\;\; \bigg |\;\;
   \begin{array}{l}
   (x,y)_{\c H}=(x,z)_{\c H}=(y,z)_{\c H}=0,\\[3mm]
   (x,x)_{\c H}=(y,y)_{\c H}=(z,z)_{\c H},
   \end{array}
\end{equation}
where $(\cdot,\cdot)_{\c H}$ is the canonical pairing
defined above.
Let $D\subset H^2(M)\times H^2(M)\times H^2(M)$ be the set
defined by the equations \eqref{_image_of_P_hyp_Equation_}.
Using Torelli theorem and Calabi-Yau, we prove the following
statement:

\hfill

\theorem\label{_image_of_P_hyp_Theorem_} 
The image of $P_{hyp}$ is Zariski dense in $D$.

\hfill

\ref{_image_of_P_hyp_Theorem_} shows that all algebraic relations
which are true for \[ (x,y,z)\in P_{hyp}(Hyp) \] are true
for all $(x,y,z)\in D$. Computing the Lie algebra $\goth a$
as in \ref{_so_5_Theorem_}, we obtain a number of relations
between $x,y,z\in H^2(M)$ which hold for all $(x,y,z)\in Im(P_{hyp})$.
Using the density argument, we obtain that these relations
are universally true. This idea leads to the proof of
\ref{_structure_alge_for_coho_hyperkahe_Theorem_}.

\hfill

The proof of \ref{_H-R_form_defi_Theorem_} is deduced from the
standard period argument and \ref{_twistor_connect_Theorem_}.
Let $\c H$ be a hyperk\"ahler structure corresponding to $I$ and $\omega$.
Clearly from the definition, the form $(\cdot,\cdot)_{\c H}$
depends only from the twistor line $\c H$, and not from the
choice of particular $I$ and $\omega$. A computation shows that
$(\cdot,\cdot)_{\c H}$  depends from $P(I)$ and not from $\omega$.
Using the fact that $C$ is all connected with twistor lines
(\ref{_twistor_connect_Theorem_}\footnote{In \cite{_main_}, we
proved a slightly weaker statement, which still suffices to prove
\ref{_H-R_form_defi_Theorem_}.}), we prove that
$(\cdot,\cdot)_{\c H}$is independent from $\c H$.


\section{Implications.}


This section contains implications of our results.

\hfill

{\bf 5.1 Mirror symmetry.} (\cite{_V:Mirror_})
Using \ref{_embe_Mum_Tate_to_so_Theorem_}
and \ref{_embe_Mum_Tate_to_so_Theorem_}, we compute
the variation of Hodge structures corresponding to the
universal VHS over the moduli space $Comp$. In
\cite{Verbitsky:Symplectic_II_}, it is proven that
for ``sufficiently generic'' deformation $W$ of a given
compact holomorphically symplectic manifold $M$, the manifold
$W$ admits no closed holomorphic curves. Therefore, using the definition
of quantum cohomology from \cite{_Kontsevich-Manin_}, we
can easily compute the quantum variation of Frobenius algebras.
Comparing these computations, we find
that Mirror Conjecture is true for holomorphically symplectic manifolds,
which are Mirror self-dual.

In proof of Mirror Symmetry, we use the fact that
tangent bundle $TM$ of a holomorphically symplectic manifold
is isomorphic to the cotangent bundle $\Omega^1(M)$ thereof. For every
Calabi-Yau manifold $M$, $\dim M =n$, the Serre's duality induces an
 isomorphism\footnote{Canonical up to a choice of a
non-degenerate section of $\Omega^n(M)$.}

\begin{equation}\label{_Yukawa_isomo_Equation_}
   H^{p}(\Omega^q(M)) \cong H^{p}(\Lambda^{n-q}(TM))
\end{equation}
beweeen cohomology of the holomorphic differential forms
and cohomology of exterrior powers of holomorphic tangent bundle.
Using the isomorphism $TM\cong \Omega^1(M)$, we interpret
the isomorphism \eqref{_Yukawa_isomo_Equation_} as a map $\eta$
from the total cohomology space $H^*(M)$ to itself. A linear-algebraic
check ensures that this map is involutive. A slightly less elementary
consideration shows that $\eta:\; H^*(M)\arrow H^*(M)$ belongs to the Lie
group $G\subset End(H^*(M))$ corresponding to the Lie algebra $\g(A)$
from \ref{_structure_alge_for_coho_hyperkahe_Theorem_}.
Clearly, Yukawa multiplication is equal to the cup-product
in cohomology twisted by $\eta$. This gives a way to
describe Yukawa product explicitely in terms of Lie algebra action.

\hfill

{\bf 5.2 Twistor paths.}

\definition 
Let $M$ be a holomorphically symplectic manifold,
$Comp$ be its moduli space, $P_0$, ... $P_n\subset Comp$ be
a sequence of twistor lines, supplied with an intersection point
$x_{i+1}\in P_i\cap P_{i+1}$ for each $i$. We say that
$\gamma= P_0, ..., P_n, x_1, ..., x_n$ is
a {\bf twistor path}. Let $I$, $I'\in Comp$.
We say that $\gamma$ is {\bf a twistor path
connecting $I$ to $I'$} if $I\in P_0$ and $I'\in P_n$.
The lines $P_i$ are called {\bf the edges},
and the points $x_i$ {\bf the vertices}
of a twistor path.

\hfill

\ref{_twistor_connect_Theorem_} proves that
every two points $I$, $I'$ in $Comp$ are connected
with a twistor path. Clearly, each twistor path
induces a diffeomorphism $\mu_\gamma:\; (M,I)\arrow (M,I')$.
We are interested in algebro-geometrical properties of this
diffeomorphism.

\hfill

For every hyperk\"ahler structure $\c H$ on $M$, let
$\g_{\c H}\subset End(H^*(M))$ be the corresponding $\goth{su}(2)$
embedded to $End(H^*(M))$. Let $H^*(M)^{\g_{\c H}}$ be the
$\g_{\c H}$-invariant part of $H^*(M)$.
Let $I\in Comp$ and $\c H$ be a hyperk\"ahler structure
which induces $I$. We say that $I$ is {\bf of general type
with respect to $\c H$} if

\[ H^*(M)^{\g_{\c H}}\cap H^*(M, \Z) =
   \oplus H^{p,p} \cap H^*(M, \Z).
\]
In \cite{Verbitsky:Symplectic_II_}, we prove that for every hyperk\"ahler
structure, all induced complex structures are of general type,
except may be a countable number thereof. Results
of \cite{_Verbitsky:Hyperholo_bundles_} and
\cite{Verbitsky:Symplectic_II_} can be compressed down to
the following statement.

\hfill

\theorem\label{_generi_implication_Theorem_} 
Let $\c H$ be a hyperk\"ahler structure on $M$ and $I$  be an
induced complex structure of general type.

(i) (\cite{Verbitsky:Symplectic_II_})
Let $N$ be a closed complex analytic subset of $(M, I)$. Then
$N$ is complex analytic with respect to $J$, for all induced
complex structures%
\footnote{In \cite{Verbitsky:Symplectic_II_},
such subsets are called {\bf trianalytic}.}%
$J$.

(ii) (\cite{_Verbitsky:Hyperholo_bundles_})
Let $Bun_I$ be the tensor category of
polystable\footnote{Polystable means direct sum of stable.
Stability is understood in the sense of Takemoto -- Mumford.}
holomorphic vector bundles of slope $0$ over $(M,I)$. For arbitrary
induced complex structure $J$, there exist a natural
injective tensor functor
$\Phi_{I\rightarrow J}:\; Bun_I\arrow Bun_J$, which is an equivalence
of $J$ is of general type with respect to $\c H$.
For $I, J, J'$ being induced complex structures and
$I$, $J$ of general type, we have

\[ \Phi_{I\rightarrow J}\circ \Phi_{J\rightarrow J'}
   = \Phi_{I\rightarrow J'}.
\]

{\bf Remark on proof of
\ref{_generi_implication_Theorem_} (ii):}
\ref{_generi_implication_Theorem_} (ii) is an implication
of the following result from \cite{_Verbitsky:Hyperholo_bundles_}.
Let $B$ be a polystable bundle on a holomorphically symplectic
Kaehler manifold $M$. We associate with the Kaehler structure
on $M$ a canonical hyperkaehler structure $\c H$
as in Calabi-Yau theorem. Assume that the first and second
Chern classes of stable summands of $B$ are
invariant under the natural action of $SU(2)$
in cohomology.
Then there exist a unique
holomorphic connection on $B$ which is holomorphic
under each of complex structures induced by $\c H$.
This lets one identify the categories
of polystable bundles for different complex structures $L$
induced by $\c H$, provided that $L$ is of general type with respect
to $\c H$.

\hfill

\definition 
Let $I$, $J\in Comp$ and $\gamma= P_0, ... P_n$ be a twistor path from
$I$ to $J$, which corresponds to the hyperk\"ahler structures
$\c H_0$, ..., $\c H_n$. We say that $\gamma$ is admissible
if $I$ is of general type with respect to $P_0$, $J$ to $P_n$,
and all vertices of $\gamma$ are of general type with respect
to the corresponding edges.

\hfill

\corollary \label{_admi_twi_impli_Corollary_} 
Let $I$, $J\in Comp$, and $\gamma$ be admissible twistor path from
$I$ to $J$.

(i) Let $\mu_\gamma:\; (M, I) \arrow (M, J)$ be the
corresponding diffeomorphism. Then, for every
complex analytic subset $N \subset (M, I)$,
$\mu_\gamma(N)$ is complex analytic with respect
to $J$, for all induced complex structures.

(ii) There exist a natural
isomorphism of tensor cetegories
\[ \Phi_{\gamma}:\; Bun_I\arrow Bun_J.\]

{\bf Proof:} Follows from \ref{_generi_implication_Theorem_}.
$\;\;\blacksquare$

\hfill

To sum it up, whenever we can connect two complex structures
by an admissible twistor path, these complex structures are
quite similar from algebro-geometrical point of view.
There is a cohomological criterion of existence of
admissible twistor path,
which is proven in the similar fashion to
\ref{_twistor_connect_Theorem_}.

\hfill

\newcommand{\Q}{{\Bbb Q}}

For $I\in Comp$, denore by $\mbox{NS}(I, \Q)$ the space
$H^{1,1}(M, I)\cap H^2(M, \Q)\subset H^2(M)$.
Let $Q\subset H^2(M, \Q)$ be a subspace of $H^2(M, \Q)$.
Let
\[ Comp_Q:= \{ I\in Comp \;\; | \;\; \mbox{NS}(I, \Q) =Q\}. \]

\hfill

\theorem\label{_admi_exi_Theorem_} 
Let $\c H$, $\c H'$
be hyperk\"ahler structures, and
$I$, $I'$ be complex structures of general type to and induced by
$\c H$, $\c H'$. Assume that $\mbox{NS}(I, \Q) = \mbox{NS}(I', \Q) =Q$,
and $I$, $I'$ lie in the same connected component of $Comp_Q$.
Then $I$, $I'$ can be connected by an admissible path.

{\bf Proof:} Follows the proof of \ref{_twistor_connect_Theorem_}.
$\;\;\blacksquare$

\hfill

For general $Q$, we have no control over the number of connected
components of $Comp_Q$ (unless global Torelli theorem is proven),
and therefore we cannot directly apply
\ref{_admi_exi_Theorem_} to obtain results from algebraic
geometry.\footnote{Exception is K3 surface, where Torelli holds.
For K3, $Comp_Q$ is connected for all $Q\subset H^2(M, \Q)$.}
However, when $Q=\emptyset$, $Comp_Q$ is clearly
connected and open in $Comp$, assuming that $Comp$ is connected,
which we assumed. On the other hand, for $I\in Comp_\emptyset$,
and every $\c H$ inducing $I$, $I$ is of general type with respect
to $\c H$ (this is essentially an implication of
\ref{_embe_Mum_Tate_to_so_Theorem_}). This proves the following corollary.

\hfill

\corollary 
Let $I$, $I'\in Comp_\emptyset$. Then $I$ can be connected to $I'$
by an admissible twistor path.

$\blacksquare$

\hfill

{\bf Remark:}
We obtain that for all $I\in Comp_\emptyset$, the closed complex
analytic subsets of $(M, I)$ have the same real analytic
structure, and categories of polystable holomorphic vector
bundles are isomorphic. There are non-trivial polystable
holomorphic vector bundles over such manifolds (tangent bundle
and its tensor powers come to mind). It is not completely clear
if manifolds $(M, I)$ with $I\in Comp_\emptyset$ have any closed
complex analytic subvarieties, except points.

\hfill

{\bf 5.3 Generalization of $(\cdot,\cdot)_{\c H}$.}
Unlike the (otherwise clearly superior)
approach used by Beauville and Bogomolov, our way of constructing
the form $(\cdot,\cdot)_{\c H}$ lends itself to an immediate
generalization.
Let $\g_0(A)$ be the grading-zero part of $\g(A)$ computed in
\ref{_Mum_Tate_is_g_0_Theorem_}, and $H^*(M)^{\g_0(A)}$
be the space of all vectors invariant under $\g_0(A)$.
Let $H^*_{\bf r}(M)$ be a subalgebra of cohomology
generated by $H^2(M)$ and $H^*(M)^{\g_0(A)}$.%
\footnote{There are only two known series of
compact hyperk\"ahler manifolds: Hilbert schemes of
Artinian sheaves on K3 surfaces, and Hilbert schemes
of Artinian sheaves on compact 2-dimensional tori, factorized by free
action of a compact torus. In both cases, the cohomology algebra
is computed by Nakajima (\cite{_Nakajima_}). It seems reasonable
to conjecture that, in either of these cases, $H^*(M)=H^*_{\bf r}(M)$.}
Let $\c H$ be a hyperk\"ahler structure on $M$. Consider the
corresponding action of $SU(2)$ on $H^*(M)$.
Let $H^i(M)= \oplus_w H^i_w(M)$ be an isotypic decomposition
of $H^i(M)$ corresponding to this action. By definition, $H^i_w(M)$ is a
direct sum of isomorphic $SU(2)$-representation of weight $w$,
where $w$, $0\leq w \leq i$ runs through the natural numbers
of the same parity as $i$. Let $(\cdot,\cdot)_{Her}$
be the Hermitian metrics on cohomology induced by the Riemannian
structure on $M$, and $(\cdot,\cdot)_{\c H}$ be the pairing which
is equal to $(-1)^{\frac{i-w}{2}} (\cdot,\cdot)_{Her}$ on $H^i_w(M)$.

\hfill

\theorem 
Consider restriction of $(\cdot,\cdot)_{\c H}$ to $H^*_{\bf r}(M)$.
This restriction $(\cdot,\cdot)_{\c H}$ is non-degenerate and
independent on $\c H$ (up to a constant multiplier).

\hfill

{\bf Proof:} For $i=2$, this statement coinsides with the statement
of \ref{_H-R_form_defi_Theorem_}. For
general $i$, the proof is essentially linear-algebraic and
identical to the proof of \cite{_main_}, Theorem 6.1.
$\;\;\blacksquare$

\hfill

{\bf Acknowledgements:}
I am grateful to my advisor David Kazhdan for warm support
and encouragement, P. Deligne and F. Bogomolov for
their suggestions and correstions, D. Kaledin and T. Pantev
for stimulating discussions. Also I owe to
Julie Lynch and IP Press for providing me with
employment.

\hfill

\end{document}